\documentclass{amsart}
\usepackage{amssymb}
\usepackage{hyperref}
\usepackage{amscd}
\newcommand{\bbc}{\mathbb{C}}\newcommand{\bbf}{\mathbb{F}}
\newcommand{\bbg}{\mathbb{G}}
\newcommand{\bbq}{\mathbb{Q}}

\newcommand{\bbz}{\mathbb{Z}}
\newcommand{\caa}{\mathcal{A}}\newcommand{\cac}{\mathcal{C}}
\newcommand{\caf}{\mathcal{F}}\newcommand{\cag}{\mathcal{G}}
\newcommand{\cah}{\mathcal{H}}\newcommand{\caj}{\mathcal{J}}
\newcommand{\cak}{\mathcal{K}}\newcommand{\call}{\mathcal{L}}

\newcommand{\cao}{\mathcal{O}}
\newcommand{\car}{\mathcal{R}}
\newcommand{\cas}{\mathcal{S}}\newcommand{\cat}{\mathcal{T}}
\newcommand{\cau}{\mathcal{U}}\newcommand{\cav}{\mathcal{V}}
\newcommand{\cax}{\mathcal{X}}\newcommand{\cay}{\mathcal{Y}}\newcommand{\caz}{\mathcal{Z}}

\newcommand{\frA}{\mathfrak{A}}

\newcommand{\frF}{\mathfrak{F}}

\newcommand{\frO}{\mathfrak{O}}
\newcommand{\frq}{\mathfrak{q}}
\newcommand{\frS}{\mathfrak{S}}\newcommand{\frX}{\mathfrak{X}}
\newcommand{\bfA}{\mathbf{A}}

\DeclareMathOperator{\aut}{Aut}
 
\DeclareMathOperator{\Div}{Div}\DeclareMathOperator{\codim}{codim}
\DeclareMathOperator{\End}{End}\DeclareMathOperator{\corank}{corank}
 \DeclareMathOperator{\et}{\text{\'et}}\DeclareMathOperator{\Frob}{Frob}
\DeclareMathOperator{\gal}{Gal}\DeclareMathOperator{\gl}{GL}

\DeclareMathOperator{\gsp}{GSp}
\DeclareMathOperator{\Ker}{Ker}
\DeclareMathOperator{\Hom}{Hom}

\DeclareMathOperator{\ns}{NS}
\newcommand{\ov}{\overline}
\DeclareMathOperator{\pic}{Pic}

\DeclareMathOperator{\cork}{cork}\DeclareMathOperator{\red}{red}
\DeclareMathOperator{\rk}{rank}
\DeclareMathOperator{\spec}{Spec}
\newfont{\cyrm}{wncysc10}

\DeclareMathOperator{\ur}{ur}

\newtheorem{theorem}[subsection]{Theorem}
\theoremstyle{plain}

\newtheorem{lemma}[subsection]{Lemma}\newtheorem{proposition}[subsection]{Proposition}
\theoremstyle{definition}\newtheorem{definition}[subsection]{Definition}
\newtheorem{remark}[subsection]{Remark}\numberwithin{equation}{subsection}
\newtheorem{hypothesis}[subsection]{Hypothesis} 
\begin{document}
\title[Selmer groups]{Selmer groups of abelian varieties in extensions of function fields}
\author{Am\'{\i}lcar Pacheco}
\address{Universidade  Federal do Rio de Janeiro, Instituto de Matem\'atica, 
Rua Guaiaquil  83, Cachambi, 20785-050 Rio de Janeiro, RJ, Brasil}
\email{amilcar@acd.ufrj.br}
\thanks{This work was partially supported by CNPq research grant 305731/2006-8.
I would like to thank Jordan Ellenberg for comments
  and suggestions on a previous version of this paper. I also thank the referee for the careful reading and suggestions which helped to improve the paper.}
\date{March 15, 2008}
\begin{abstract}
Let $k$ be a field of characteristic $q$, $\cac$ a smooth geometrically connected
curve defined over $k$ with function field $K:=k(\cac)$. Let $A/K$ be
a non constant abelian variety defined over $K$ of dimension $d$. We
assume that $q=0$ or $>2d+1$. Let $p\ne q$
be a prime number and $\cac'\to\cac$ a finite geometrically \textsc{Galois} and
\'etale cover defined over $k$ with function field $K':=k(\cac')$. 
Let $(\tau',B')$ be the
$K'/k$-trace of $A/K$. We give an upper bound for the $\bbz_p$-corank
of the \textsc{Selmer} group $\text{Sel}_p(A\times_KK')$, defined in terms of
the $p$-descent map. As a consequence, we get
an upper bound for the $\bbz$-rank of the \textsc{Lang-N\'eron} group
$A(K')/\tau'B'(k)$. In the case of a geometric tower of 
curves whose \textsc{Galois} group is isomorphic to $\bbz_p$, we give sufficient 
conditions for the \textsc{Lang-N\'eron} group of $A$ to be uniformly bounded along the
tower. 
\end{abstract}
\maketitle
\section{Introduction}
Let $\cac$ be a smooth
geometrically connected (not necessarily complete) curve defined over a field 
$k$ of characteristic $q\ge0$. Let $\cax$ be a regular 
compactification of $\cac$. Denote by $K:=k(\cax)=k(\cac)$ the
function field of $\cax$ (or equivalently of $\cac$) and let $g$ be
the genus of $\cax$. Let $k^s$ be a separable
closure of $k$, 
$\cak:=k^s(\cax)=k^s(\cac)$ and $\cak^s$ a separable closure of $\cak$. 

Let $A/K$ be a non constant 
abelian variety of dimension $d$ defined over $K$. A model for $A/K$
over $k$ consists of a smooth geometrically connected projective variety $\caa$
defined over $k$ and a proper flat morphism $\phi:\caa\to\cac$ also defined
over $k$ whose generic fiber is $A/K$. Denote by $\cau$ the smooth locus of
$\phi$ and $\caa_{\cau}:=\phi^{-1}(\cau)$. The morphism $\phi$ induces
an abelian scheme $\phi_{\cau}:\caa_{\cau}\to\cau$ having still $A/K$
as its generic fiber.

Let $(\tau,B)$ be the $K/k$-trace 
of $A$  and $d_0:=\dim(B)$.  
A theorem of \textsc{Lang} and \textsc{N\'eron} \cite[theorem 2,
chapter 6]{la} states the group $A(\cak)/\tau B(k^s)$ is finitely
generated. \emph{A fortiori}, the same holds for $A(K)/\tau
B(k)$. 

Let $\frF_A$ be the conductor divisor of $A/K$ in $\cax$ and $f_A$ its degree. 
\textsc{Ogg} obtained in \cite[theorem 2]{ogg} the following geometric 
upper bound for the rank of $A(\cak)/\tau B(k^s)$ 
\begin{equation}\label{bdogg}
\rk\left(\frac{A(\cak)}{\tau B(k^s)}\right)\le2d(2g-2)+f_A+4d_0.
\end{equation} 
In particular, this is also an upper bound for the rank of $A(K)/\tau
B(k)$. 

Let $\psi:\cac'\to\cac$ be a finite morphism defined over $k$ which is
also a geometrically \textsc{Galois} cover. Let $\cag:=\aut_{k^s}(\cac'/\cac)$
be its \textsc{Galois} group. The absolute \textsc{Galois} group $G_k:=\gal(k^s/k)$ of
$k$ acts on $\cag$ and we denote by $\frO_{G_k}(\cag)$ the set of
orbits of this action. 

Let $K':=k(\cac')$, $g'$ the genus of $\cac'$ and $(\tau',B')$ the
$K'/k$-trace of $A$. 
If $k$ is a number field, $\cac$ is complete
and $\psi$ is geometrically abelian, we obtained in \cite[theorem
1.3]{pac} the following upper bound
\begin{equation}\label{bdpac}
\rk\left(\frac{A(K')}{\tau'B'(k)}\right)\le\frac{\#\frO_{G_k}(\cag)}{\#\cag}(d(2d+1)(2g'-2))+
\#\frO_{G_k}(\cag)(2df_A).
\end{equation}
This result was proved under the assumptions that the \textsc{Tate}'s
conjecture for divisors holds for $\caa/k$ and that the monodromy
representations with respect to the first and second \'etale
cohomology groups of $A/K$ with coefficients in $\bbq_{\ell}$ are
irreducible. For an explanation of these two hypotheses see
\emph{loc. cit.} sections 2 and 3. In the particular case where $A$ is
the Jacobian variety of a curve and $\psi$ is \'etale, then the upper
bound in 
(\ref{bdpac}) is indeed better than that in (\ref{bdogg}) (\emph{loc. cit.}
paragraph after remark 1.8, in particular formula (1.8)).

In this paper we obtain an improvement the bound (\ref{bdpac}) in theorem
\ref{mainthm1} without
appeal to the previous hypotheses and of the same order of magnitude
as that of \cite[(1.8)]{pac}, where it was only obtained in the case where $A$
was a Jacobian variety. This extends a previous result of
\textsc{Ellenberg} in \cite[theorem 2.8]{ell} from elliptic to abelian
fibrations. In order to do this we start with descent maps in \S2,
then in \S3 we treat \textsc{Selmer} groups and prove our first result. In \S4 we give conditions
for the \textsc{Lang-N\'eron} groups of $A$ to be uniformly bounded
along a geometric tower of curves over $k$ whose \textsc{Galois} group
is isomorphic to $\bbz_p$. Our
theorem \ref{fginfty} extends \cite[theorem 4.4]{ell}. In the final
section (\S5) we discuss a sufficient condition to obtain an example (related to Jacobian fibrations) where
the hypotheses of theorem \ref{fginfty} are met.

\section{Descent maps}\label{desc}
Let $\cac_s:=\cac\times_{\spec(k)}\spec(k^s)$, 
$\cau_s:=\cau\times_{\spec(k)}\spec(k^s)$ and $\eta_s$ the geometric
generic point of $\cau$. Let $\pi_1(\cau_s,\eta_s)$ be the \emph{algebraic
fundamental group} of $\cau_s$ with respect to $\eta_s$. Let $\cak^{\ur}$ the maximal
\textsc{Galois} subextension of $\cak^s/\cak$ that is unramified over 
$\cau_s$. Then $\pi_1(\cau_s,\eta_s)\cong\gal(\cak^{\ur}/\cak)$. Let 
$\pi_1^t(\cau_s,\eta_s)$ be the \emph{tame algebraic fundamental group}  
of $\cau_s$ with respect to $\eta_s$. Let $\cak^t$ be the maximal
\textsc{Galois} subextension of $\cak^s/\cak$ which is unramified over $\cau_s$
and at most tamely ramified over $\cac_s\setminus\cau_s$. Then
$\pi_1^t(\cau_s,\eta_s)\cong\gal(\cak^t/\cak)$. 

Let $p\ne q:=\text{char}(k)$ a prime number. For each integer $n\ge1$,
let $A[p^n]$ be 
the subgroup of $p^n$-torsion points of $A$ and 
$A[p^{\infty}]:=\bigcup_{n\ge1}A[p^n]=\varinjlim_nA[p^n]$ the $p$-divisible subgroup of
$A$. Let $\cak(A[p^n])$ be the
subfield of $\cak^s$ generated over $\cak$ by the coordinates of the
points in $A[p^n]$. By the
\textsc{N\'eron-Ogg-Shafarevich} criterion \cite[theorem
1]{seta} for every $n\ge1$ the group $A[p^n]$ is unramified at every $v\in\cau_s$,
therefore $\cak(A[p^n])\subset\cak^{\ur}$. 

For each $v\in\cac_s$, let $\cak_v$ be the completion of $\cak$ at
$v$. Let $v^s$ be the unique extension of $v$ to $\cak^s$ and
$\cak^s_{v^s}$ the completion of $\cak^s$ at $v^s$. The inertia group
$I_v$ at $v$ equals the decomposition group 
$D_v$ at $v$, since $k^s$ is separably closed. Moreover, $D_v$ is
isomorphic to $\gal(\cak^s_{v^s}/\cak_v)$.  

\begin{lemma}\label{lemur}
For every integer $n\ge1$, we have a short exact sequence of groups 
\begin{equation}\label{urseqtr0}
0\to A[p^n]\to A(\cak^{\ur})\overset{\times p^n}\longrightarrow A(\cak^{\ur})\to0.
\end{equation}
\end{lemma}

\begin{proof}
It suffices to prove the surjectivity of the multiplication by $p^n$
map. Let $x\in A(\cak^{\ur})$ and $y\in A(\cak^s)$ such that
$p^ny=x$. Let $v\in\cau_s$ and 
$\sigma\in I_v$. Then $\sigma(p^ny)=p^n\sigma(y)=\sigma(x)=x=p^ny$,
thus $\sigma(y)-y\in A[p^n]$. 

Let $\call:=\cak_v((1/p^n)x)$ be the \textsc{Galois} subextension of $\cak_{v^s}^s/\cak_v$
generated by all the solutions $y'\in A(\cak_{v^s}^s)$ of $p^ny'=x$ over $\cak_v$. 
Let $w$ be the extension of $v$
to $\call$ and $\cao_w$ its valuation ring. 

Let $\bfA_w$ be the \textsc{N\'eron} minimal 
model of $A$ at $w$ and $\tilde{\bfA}_w$ its special fiber defined over $k^s$. By the
properties of the \textsc{N\'eron} model we have
$A(\call)\cong\bfA_w(\cao_w)$. Thus we have a reduction map
$\red_w:A(\call)\to\tilde{\bfA}_w(k^s)$ and we denote
$\tilde{y}:=\red_w(y)$. The restriction $\sigma_w$ of $\sigma$ to
$\call$ is an element of the inertia group $I(w|v)$ of $w$ over
$v$. Notice that $\gal(\call/\cak_v)$ is equal to the decomposition
group $D(w|v)$ of $w$ over $v$. In fact, this group coincides with
$I(w|v)$, because $\kappa_w=\kappa_v=k^s$, where $\kappa_w$,
resp. $\kappa_v$ denotes the residue field of $w$, resp. $v$. 
Since the reduction $\tilde{\sigma}_w$ of $\sigma_w$ modulo $w$
is trivial and $\red_w$ commutes with the action of
$\gal(\call/\cak_v)$, then
$\widetilde{\sigma_w(y)}=\tilde{\sigma}_w(\tilde{y})=\tilde{y}$. 
But by \cite[lemma
2 (b), chapter 7, \S7.3]{bolura} the multiplication by $p^n$ is
\'etale in $\bfA_w$. Thus $\red_w$ is injective in
$A(\call)\cap A[p^n]$, whence $\sigma(y)=\sigma_w(y)=y$ for every $v\in\cau_s$ and
$\sigma\in I_v$, \emph{i.e.}, $y\in A(\cak^{\ur})$. 
\end{proof}

\begin{remark}
Note that the exact sequence (\ref{urseqtr0}) remains exact after
passing to the quotient by $\tau B(k^s)$, thus we obtain for every
$n\ge1$ the exact sequence of groups 
\begin{equation}\label{urseq}
0\to\left(\frac{A(\cak^{\ur})}{\tau
    B(k^s)}\right)[p^n]\to\frac{A(\cak^{\ur})}{\tau
    B(k^s)}\overset{\times
    p^n}\longrightarrow\frac{A(\cak^{\ur})}{\tau B(k^s)}\to0,
\end{equation} 
where the first group of (\ref{urseq}) denotes the $p^n$-torsion subgroup of
$A(\cak^{\ur})/\tau B(k^s)$. Moreover, the long exact \textsc{Galois} cohomology
sequence derived from (\ref{urseq}) produces $p^n$-descent map 
$$
\delta_{p^n}:\frac{A(\cak)/\tau B(k^s)}{p^n(A(\cak)/\tau B(k^s))} 
\hookrightarrow
H^1\left(\pi_1(\cau_s,\eta_s),\left(\frac{A(\cak^{\ur})}{\tau B(k^s)}\right)[p^n]\right)
$$
for every $n\ge1$. Since $(A(\cak)/\tau B(k^s))/p^n(A(\cak)/\tau B(k^s))\cong
(A(\cak)/\tau B(k^s))\otimes_{\bbz}\bbz/p^n\bbz$, taking the injective limit we have
a $p^{\infty}$-descent map
$$
\delta_{p^{\infty}}:\frac{A(\cak)}{\tau B(k^s)}\otimes_{\bbz}\bbq_p/\bbz_p\hookrightarrow
H^1\left(\pi_1(\cau_s,\eta_s),\left(\frac{A(\cak^{\ur})}{\tau B(k^s)}\right)[p^{\infty}]\right). 
$$
\end{remark}

\begin{remark}\label{remsst}
For each $v\in\cac_s\setminus\cau_s$, let $R_v\subset I_v$
be the first ramification group at $v$. Then $R_v$ is the \textsc{Sylow} $q$-subgroup of
$I_v$. \textsc{Grothendieck} showed in \cite[exp. 9, (4.6.3)]{sga7} that if
$A$ has semi-stable reduction at some $v\in\cac_s$, then the Swan
conductor $\delta_{p,v}(A)$ of $A$ at $v$ with respect to a prime number
$p\ne q$ vanishes. This is equivalent to $R_v$ acting trivially on
$A[p]$. Since for every $v\in\cau_s$, the group $I_v$ acts trivially 
on $A[p]$, we conclude that $\delta_{p,v}(A)=0$ for every $v\in\cac_s$
if and only if 
$\cak(A[p])/\cak$ is unramified at every $v\in\cau_s$ and at most tamely
ramified at every $v\in\cac_s\setminus\cau_s$, \emph{i.e.},
$\cak(A[p])\subset\cak^t$. 

If $q=0$ or $q>2d+1$, \textsc{Ogg} observed in 
\cite[remark 1, p. 211]{ogg} there exists a prime number $p\ne q$
such that $\cak(A[p])/\cak$ is tamely ramified. Thus
$\delta_{p,v}(A)=0$ for every $v\in\cac_s$. 
But \textsc{Grothendieck} also showed in \cite[exp. 9, corollaire
4.6]{sga7} that $\delta_{p,v}(A)$ does not depend on $p$. In
particular, $\cak(A[p])/\cak$ is tamely ramified for every prime
number $p\ne q$. Furthermore, \textsc{Deschamps} proved in
\cite[corollaire 5.18]{mds} that
$A$ acquires semi-stable reduction over $\cak(A[p])$.
\end{remark}

\begin{hypothesis}
Summing-up, we assume from now on that $q=0$ or $q>2d+1$, so that $A$
acquires semi-stable reduction over a finite tamely ramified extension
of $\cak$, namely $\cak(A[p])$ for any prime number $p\ne q$. We also
fix from now on a prime number $p\ne q$.
\end{hypothesis}

\begin{lemma}\label{lemsst}
For every $n\ge1$ and $v\in\cac_s\setminus\cau_s$, the ramification
group $R_v$ acts trivially on $A[p^n]$. Moreover, we have an exact
sequence of groups 
\begin{equation}\label{tseqtr0} 
0\to A[p^n]\to A(\cak^t)\overset{\times p^n}\longrightarrow
A(\cak^t)\to 0.
\end{equation} 
\end{lemma}

\begin{proof}
We proceed similarly to the proof of lemma \ref{lemur}. We show
the first statement by induction. Let $y\in A[p^2]$, $x:=py\in A[p]$,
$v\in\cac_s\setminus\cau_s$ and $\sigma\in R_v$. Then
$\sigma(py)=p\sigma(y)=\sigma(x)=x=py$, \emph{i.e.}, $\sigma(y)-y\in
A[p]$. 

Let $\call:=\cak_v((1/p)x)$ be the \textsc{Galois}
subextension of $\cak_{v^s}^s/\cak_v$ generated by the solutions
$y'\in A(\cak_{v^s}^s)$ 
of $py'=x$ over $\cak_v$. Let $w$ be the extension of $v$ to $\call$. 
Let $\tilde{\bfA}_w^0$ be the connected component of the special fiber
$\tilde{\bfA}_w$ of the \textsc{N\'eron} minimal model $\bfA_w$ of $A$ over
$\cao_w$. Again the reduction $\tilde{\sigma}_w$ of the restriction
$\sigma_w$ of $\sigma$ to $\call$ modulo $w$ is trivial, thus
$\widetilde{\sigma_w(y)}=\tilde{\sigma}_w(\tilde{y})=\tilde{y}$. Once more by \cite[lemma 2
(b), chapter 7, \S7.3]{bolura} 
the map $\red_w:A(\call)\to\tilde{\bfA}_w^0(k^s)$ is injective in
$A(\call)\cap A[p]\subset A(\cak^t)$, 
whence $\sigma(y)=\sigma_w(y)=y$ for every $v\in\cac_s\setminus\cau_s$ and
$\sigma\in R_v$. We have already observed in the proof of lemma
\ref{lemur} that $\sigma(y)=y$ for every $v\in\cau_s$ and $\sigma\in
I_v$. Consequently $y\in A(\cak^t)$. Now the first claim follows by induction.

For the proof of the surjectivity of the multiplication by $p^n$ map
follows exactly the proof of lemma \ref{lemur}, once we know that for
every $v\in\cac_s\setminus\cau_s$ and $\sigma\in R_v$, $\sigma$ acts trivially on
$A[p^n]$; and similarly for every $v\in\cau_s$ and $\sigma\in I_v$,
$\sigma$ acts trivially on $A[p^n]$. 
\end{proof}

\begin{remark}
As before, the exact sequence (\ref{tseqtr0}) remains also exact after
taking the quotient by $\tau B(k^s)$ leading to the exact sequence
\begin{equation}\label{tseq}
0\to\left(\frac{A(\cak^t)}{\tau
    B(k^s)}\right)[p^n]\to\frac{A(\cak^t)}{\tau B(k^s)}\overset{\times
    p^n}\longrightarrow\frac{A(\cak^t)}{\tau B(k^s)}\to0.
\end{equation}
The long exact \textsc{Galois} cohomology sequence derived from (\ref{tseq})
shows that $\delta_{p^n}$ factors through the tame $p^n$-descent map
$$
\delta_{p^n}^t:\frac{A(\cak)/\tau B(k^s)}{p^n(A(\cak)/\tau
  B(k^s))}\hookrightarrow
H^1\left(\pi_1^t(\cau_s,\eta_s),\left(\frac{A(\cak^t)}{\tau B(k^s)}\right)[p^n]\right).
$$
Taking the injective limit we obtain a tame $p^{\infty}$-descent map 
$$
\delta_{p^{\infty}}^t:\frac{A(\cak)}{\tau
  B(k^s)}\otimes_{\bbz}\bbq_p/\bbz_p
\hookrightarrow
H^1\left(\pi_1^t(\cau_s,\eta_s),\left(\frac{A(\cak^t)}{\tau B(k^s)}\right)[p^{\infty}]\right).
$$
\end{remark}

We end this section with following two observations.

\begin{lemma}\label{leminvpdiv}
The subgroup of the elements of $(A(\cak^t)/\tau B(k^s))[p^{\infty}]$ which are fixed under
the action of $\pi_1^t(\cau_s,\eta_s)$ is finite.
\end{lemma}

\begin{proof}
Note that the set of elements of $(A(\cak^t)/\tau B(k^s))[p^{\infty}]$
fixed by $\pi_1^t(\cau_s,\eta_s)\cong\gal(\cak^t/\cak)$ is simply
$(A(\cak)/\tau B(k^s))[p^{\infty}]$ which is finite by the
\textsc{Lang-N\'eron} theorem. 
\end{proof}

\begin{remark}
Observe that for every $n\ge1$ we have 
\begin{equation}\label{pntor}
\left(\frac{A(\cak^t)}{\tau B(k^s)}\right)[p^n]=\frac{A[p^n]+\tau
  B(k^s)}{\tau B(k^s)}\cong\frac{A[p^n]}{\tau B[p^n]}. 
  \end{equation}
Denote $A[p^{\infty}]/\tau B[p^{\infty}]:=\varinjlim_nA[p^n]/\tau
B[p^n]$. 
As a consequence we rewrite the tame $p^{\infty}$-descent map as 
\begin{equation}\label{tamepinftydesc} 
\delta_{p^{\infty}}^t:\frac{A(\cak)}{\tau
  B(k^s)}\otimes_{\bbz}\bbq_p/\bbz_p\hookrightarrow
H^1\left(\pi_1^t(\cau_s,\eta_s),
\frac{A[p^{\infty}]}{\tau B[p^{\infty}]}\right). 
\end{equation} 
\end{remark} 

\section{Selmer groups}\label{selgp}

\begin{definition}
Let $j:\eta_s\hookrightarrow\cac_s$ be the inclusion map. For each discrete
$p$-primary torsion \'etale sheaf $F$ in $\eta_s$, let
$\caf:=j_*F$. The \textsc{Selmer} group of $\cac$ with respect to $F$, denoted
by $S(\cac,F)$, is defined as the first \'etale cohomology group
$H^1_{\et}(\cac_s,\caf)$ of $\cac_s$ with coefficients in $\caf$. Consider the
$p$-primary \'etale sheaf $F_{p^{\infty}}:=A[p^{\infty}]/\tau B[p^{\infty}]$ in $\eta_s$.
Let $\caf_{p^{\infty}}:=j_*(F_{p^{\infty}})$. In the abstract, we
denoted $S(\cac,F_{p^{\infty}})$ by $\text{Sel}_p(A/K)$. 
\end{definition}

\begin{remark}Using some results on \'etale cohomology (\emph{cf.} \cite[III.1.25, V.2.17]{mil}), the Selmer group $S(\cac,F)$ can be alternatively described in terms of  \textsc{Galois} cohomology by 
\begin{equation}\label{seldef}S(\cac,F)=\Ker\left(H^1(\pi_1^t(\cau_s,\eta_s),F)\to\bigoplus_{v\in\cac_s\setminus\cau_s}H^1(\gal(\cak_v^t/\cak_v),F)\right),\end{equation}
where $\cak_v^t/\cak_v$ is the maximal \textsc{Galois} subextension of $\cak_{v^s}^s/\cak_v$ tamely ramified over $\cak_v$.
In the case where $F=F_{p^{\infty}}$, 
this definition also agrees with the
classical one for $S(\cac,F_{p^{\infty}})$. First observe that the local tame 
$p^{\infty}$-descent map 
$$
\delta_{p^{\infty},v}^t:\frac{A(\cak_v)}{\tau B(k^s)}
\otimes_{\bbz}\bbq_p/\bbz_p\hookrightarrow
H^1(\gal(\cak_v^t/\cak_v),F_{p^{\infty}}) 
$$ 
has trivial image. Indeed, by a
result of \textsc{Mattuck} \cite{mat}, $A(\cak_v)$ is
isomorphic to $\cao_v^d$ (as an additive group) times a finite
group, where $\cao_v$ denotes the ring of integers of $\cak_v$.
Since $\text{char}(\cak_v)=q\ne p$, then
$(A(\cak_v)/\tau B(k^s))\otimes_{\bbz}\bbq_p/\bbz_p=0$. 
\end{remark}

Before we proceed, let us recall the definition of the conductor
$\frF_A$ of the abelian variety $A/K$ of dimension $d$.
Let $T_{p}(A)$ the $p$-adic \textsc{Tate} 
module of $A$ and $V_{p}(A):=T_{p}(A)\otimes_{\bbz_{p}}\bbq_{p}$. 
Let $\cax_s:=\cax\times_{\spec(k)}\spec(k^s)$. 

\begin{definition}\label{defcond}
The multiplicity $e_v$ of the conductor $\frF_A$ of $A$ at $v\in\cax_s$ is
defined as $\codim(V_{p}(A)^{I_v})$, where
$V_{p}(A)^{I_v}$
denotes the subspace of $V_{p}(A)$
which is fixed by the action of the inertia group $I_v$. 
The multiplicity $e_v$ can be described in terms of
the type of the reduction of $A$ at $v$ as follows. Let
$\bfA_v\to\spec(\cao_v)$ be the \textsc{N\'eron} minimal model of $A$ at
$v\in\cax$, $\tilde{\bfA}_v\to\spec(k^s)$
its special fiber and $\tilde{\bfA}_v^0$ the connected component of
$\tilde{\bfA}_v$. Then $\tilde{\bfA}_v^0$
is an extension of an abelian variety $\frA_v$ by a linear algebraic
group $L_v$ both defined over $k^s$. Let $a_v:=\dim(\frA_v)$. The
linear algebraic group $L_v$ is equal to a product
$\bbg_{m}^{t_v}\times \bbg_{a}^{u_v}$. 
The non-negative integers $a_v$, $t_v$ and
$u_v$ are called the \emph{abelian}, \emph{reductive} and \emph{unipotent} ranks of
$\tilde{\bfA}_v^0$, resp.. Since we are supposing that $k$ has characteristic either 0 or greater than $2d+1$, we have no contribution from the \textsc{Swan} conductor, in this circumstance the multiplicity $e_v$ of $\frF_A$ at $v$ equals $t_v+2u_v$. Recall also that $d=a_v+t_v+u_v$.
\end{definition}

Let $M$ be a $\bbz_p$-module and $M^*:=\Hom(M,\bbq_p/\bbz_p)$ its
\textsc{Pontryagin} dual module. If $M^*$ is finitely generated, then we say
that $M$ is \emph{cofinitely generated}. In this case the rank of $M^*$ is
called the corank of $M$ denoted by
$\cork_{\bbz_p}(M)$. We say that $M$ is a
$\bbz_p$-\emph{cofree} module, if $M^*$ is a $\bbz_p$-free module.

\begin{lemma}\label{corka}
Let
$$
\cah(\cac,F_{p^{\infty}}):
=\bigoplus_{v\in\cac_s\setminus\cau_s}H^1(\gal(\cak_v^t/\cak_v),F_{p^{\infty}}).
$$
Then $\cah(\cac,F_{p^{\infty}})$ is a
  $\bbz_p$-cofree module of corank equal to
$$
\sum_{v\in\cac_s\setminus\cau_s}(2a_v+t_v)-\#(\cac_s\setminus\cau_s)\cdot(2d_0).
$$
\end{lemma}

\begin{proof}
Note that $\gal(\cak_v^t/\cak_v)\cong I_v/R_v$ is procyclic. 
It follows from \cite[1.6.13 (i), p. 69]{neucoh} that  
$$
H^1(\gal(\cak_v^t/\cak_v),F_{p^{\infty}})
\cong (F_{p^{\infty}})_{\gal(\cak_v^t/\cak_v)}, 
$$ 
where $(F_{p^{\infty}})_{\gal(\cak_v^t/\cak_v)}$ denotes the subspace
of $F_{p^{\infty}}$ of the elements which are coinvariant with respect
to the action of $\gal(\cak_v^t/\cak_v)$. 
But, by lemma
\ref{lemsst}, $R_v$ acts as identity on $A[p^{\infty}]$, therefore 
$$
(F_{p^{\infty}})_{\gal(\cak_v^t/\cak_v)}\cong(F_{p^{\infty}})_{I_v}
\cong(A[p^{\infty}])_{I_v}/\tau B[p^{\infty}]. 
$$ 
The latter is $\bbz_p$-cofree with $\bbz_p$-corank equal to
$2d-e_v-2d_0=2a_v+t_v-2d_0$ and this implies the lemma. 
\end{proof}

\begin{definition}
Let $\frF_{A,\cac}:=\sum_{v\in\cac}e_vv$ be the restriction of the
conductor $\frF_A$ of $A$ to $\cac$ and
$f_{A,\cac}:=\deg(\frF_{A,\cac})$. 
\end{definition}

\begin{proposition}\label{corkprop}(cf. \cite[proposition 2.5]{ell})
\begin{enumerate}
\item $H^1(\pi_1^t(\cau_s,\eta_s),F_{p^{\infty}})$ is a cofree $\bbz_p$-module of corank
  equal to
$$
(2d-2d_0)(2g-2+\#(\cax_s\setminus\cac_s))+f_{A,\cac}+\sum_{v\in\cac_s\setminus\cau_s}(2a_v+t_v)
-\#(\cac_s\setminus\cau_s)\cdot(2d_0).
$$
\item $S(\cac,F_{p^{\infty}})$ is a $\bbz_p$-module of corank
  equal to
$$
(2d-2d_0)(2g-2+\#(\cax_s\setminus\cac_s))+f_{A,\cac}.
$$
\end{enumerate}
\end{proposition}

\begin{proof}
The proof is similar to \cite[proposition 2.5]{ell}, replacing
\cite[remark 2.4]{ell} by lemma \ref{leminvpdiv}, taking into
account that \cite[chapter V, (2.18)]{mil} yields that
$H^1(\pi_1^t(\cau_s,\eta_s),F_{p^{\infty}})$ has $\bbz_p$-corank equal
to
\begin{align*}
&(2d-2d_0)(2g-2+\#(\cax_s\setminus\cau_s))\\
&=(2d-2d_0)(2g-2+\#(\cax_s\setminus\cac_s))
+\sum_{v\in\cac_s\setminus\cau_s}(2a_v+t_v)+f_{A,\cac}
-\#(\cac_s\setminus\cau_s)\cdot(2d_0). 
\end{align*}
and using lemma \ref{corka}. 
\end{proof}

The framework of lemma \ref{corka} and proposition \ref{corkprop}
allows one to immediately extend \cite[theorem 2.8]{ell}
to abelian varieties as follows. 

\begin{definition}\cite[definition 2.6]{ell}
Let $G$ be a finite group, $H$ a subgroup of $\aut(G)$ and
$\Gamma:=G\rtimes H$ the semi-direct product of $G$ and $H$.
Denote by $V_{\Gamma}$, resp. $V_G$, the real vector
space spanned by the irreducible complex-valued characters of $\Gamma$,
resp. $G$. An element $v\in V_{\Gamma}$, resp. $V_G$, is non-negative if
$\langle v,\psi\rangle\ge0$ for every character $\psi$ of an
irreducible representation of
$\Gamma$, resp. $G$. Denote by $[\Gamma/H]\in V_{\Gamma}$ the coset character of $\Gamma$
with respect to $H$ and by $[G/1]\in V_G$ the regular character of
$G$. Let $\epsilon(G,H)$ be the maximum of $\langle v,[\Gamma/H]\rangle$ over
all $v\in V_{\Gamma}$ such that
\begin{enumerate}
\item $v\ge0$.
\item $[G/1]-r(v)\ge0$, where $r:V_{\Gamma}\to V_G$ is the restriction map.
\end{enumerate}
This number is well-defined, since the region of $V_{\Gamma}$ defined
by the two previous conditions is a compact polytope.
\end{definition}

Let $\psi:\cac'\to\cac$ a finite morphism defined over $k$ which is
geometrically \textsc{Galois}. 
Let $l/k$ be the smallest finite \textsc{Galois} extension 
of $k$ such that all elements of $\cag=\aut_{\ov{k}}(\cac'/\cac)$ 
are defined over $l$. Hence, $G_k$ acts on $\cag$ via the finite
quotient $H:=\gal(l/k)$. 

\begin{theorem}\label{mainthm}(cf. \cite[theorem 2.8]{ell})
Let $\cac/k$ be a smooth geometrically connected curve of genus $g$ defined over a field $k$ with function field $K=k(\cac)$. Let $A/K$ be a non constant abelian variety of dimension $d$. Suppose that $k$ has characteristic either 0 or greater than $2d+1$. Let $d_0$ be the dimension of the $K/k$-trace of $A$ and $\cax$ a regular compactification of $\cac$. Denote by $k^s$ a separable closure of $k$, and $\cac_s:=\cac\times_{\spec(k)}\spec(k^s)$ and $\cax_s:=\cax\times_{\spec(k)}\spec(k^s)$. Let $f_{A,\cac}$ be the degree of the conductor of $A$ with respect to $\cac$. 
Let $\psi:\cac'\to\cac$ be a finite morphism defined over $k$ which is 
geometrically \textsc{Galois} and \'etale with automorphism group $\cag=\aut_{k^s}(\cac'/\cac)$.
Then
\begin{equation}\label{bdell}
\rk\left(\frac{A(K')}{\tau'B'(k)}\right)\le 
\epsilon(\cag,H)((2d-2d_0)(2g-2+\#(\cax_s\setminus\cac_s))+f_{A,\cac}).
\end{equation}
\end{theorem}

\begin{proof}
The proof is the same as in \cite[theorem
2.8]{ell} replacing \cite[lemma 2.9]{ell} by lemma \ref{qpsp}. 
\end{proof}

Let $\caa':=\caa\times_{\cac}\cac'$ and $\phi':\caa'\to\cac'$ the
morphism obtained from $\phi$ by changing the base from $\cac$ to
$\cac'$. Denote by $\cau'$ the smooth locus of $\phi'$ and let
$\cau'_s:=\cau'\times_{\spec(k)}\spec(k^s)$. Let $\eta'_s$ be a
geometric generic point of $\cau'$ and $\pi_1^t(\cau'_s,\eta'_s)$ the
tame algebraic fundamental group of $\cau_s'$ with respect to
$\eta'_s$. Let $K':=k(\cac')$ and denote by $(\tau',B')$ the $K'/k$-trace of $A$. Let $d'_0:=\dim(B')$. Let $F_{p^{\infty}}':=A[p^{\infty}]/\tau'B'[p^{\infty}]$. 

\begin{lemma}\label{qpsp} (cf. \cite[lemma 2.9]{ell})
With the same hypotheses of theorem \ref{mainthm},
$$
\Hom(S(\cac',F_{p^{\infty}}'),\bbq_p/\bbz_p)\otimes_{\bbz_p}\bbq_p
$$
is a free $\bbq_p[\cag]$-module of rank at most 
$$
(2d-2d_0)(2g-2+\#(\cax_s\setminus\cac_s))+f_{A,\cac}.
$$
\end{lemma}

\begin{proof}
For any discrete cofinitely generated $\bbz_p[\cag]$-module $M$ we
associate the \linebreak finitely generated $\bbq_p[\cag]$-module
$W(M):=\Hom(M,\bbq_p/\bbz_p)\otimes_{\bbz_p}\bbq_p$. 

Let $\call:=k^s(\cac')$, $w\in\cac'_s:=\cac'\times_{\spec(k)}\spec(k^s)$ and $\call_w$ the completion of
$\call$ at $w$. Let $\call^s$ be a separable closure of $\call$,
$w^s$ the extension of $w$ to $\call^s$ and $\call^s_{w^s}$ the
completion of $\call^s$ at $w^s$. Let $\call_w^t/\call_w$ be the
maximal \textsc{Galois} subextension of $\call_{w^s}^s/\call_w$ which is 
tamely ramified over $\call_w$.  

As before, we consider the following $\bbz_p[\cag]$-module
$$
\cah(\cac',F'_{p^{\infty}}):=\bigoplus_{w\in\cac'_s\setminus\cau'_s}
H^1(\gal(\call_w^t/\call_w),F'_{p^{\infty}}).
$$
It follows from the proof of lemma \ref{corka} that for each
$w\in\cac'_s\setminus\cau_s'$ the $\bbz_p$-module
$H^1(\gal(\call_w^t/\call),
F'_{p^{\infty}})$ has $\bbz_p$-corank
$2a_w+t_w-2d_0'$. Since $\cac'\to\cac$ is
geometrically \'etale, the morphism
$\bfA_v\times_{\spec(\cao_v)}\spec(\cao_w)\to\bfA_w$ of change of base
of \textsc{N\'eron} models 
is an isomorphism (\emph{cf.} \cite[chapter 7,
theorem 1, p. 176]{bolura}). \emph{A fortiori}, $a_w=a_v$, $t_w=t_v$
and $u_w=u_v$ for every $w\in\cac'_s$ over $v\in\cac_s$. As a
consequence, $W(\cah(\cac',F'_{p^{\infty}}))$ is a
free $\bbq_p[\cag]$-module of rank equal to 
\begin{equation}\label{weq}
\sum_{v\in\cac_s\setminus\cau_s}(2a_v+t_v)-\#(\cac_s\setminus\cau_s)\cdot(2d_0').
\end{equation} 
Let $j':\eta'_s\hookrightarrow\cac_s'$ be the inclusion map and
$\caf_{p^{\infty}}':=j'_*(F_{p^{\infty}}')$. 
As in \cite[proof of proposition 2.5]{ell} the finiteness of
$H^2_{\et}(\cac'_s,\caf'_{p^{\infty}})$ implies the equality
\begin{equation}\label{grothgp2}
[W(H^1(\pi_1^t(\cau'_s,\eta'_s),F'_{p^{\infty}}))]=
[W(S(\cac',F'_{p^{\infty}}))]+
[W(\cah(\cac',F'_{p^{\infty}}))]
\end{equation} 
in the \textsc{Grothendieck} group of the category of
$\bbq_p[\cag]$-modules. It follows from \textsc{Shapiro}'s lemma
\cite[proposition 10, p. I-12]{sercoh} that 
$$
H^1(\pi_1^t(\cau'_s,\eta'_s),F'_{p^{\infty}})=
H^1(\pi_1^t(\cau_s,\eta_s),F'_{p^{\infty}}\otimes_{\bbz}\bbz[\cag]). 
$$
By \cite[chapter V, remark 2.19]{mil}, for every
$\pi_1^t(\cau_s,\eta_s)$-module $M$ we have an identity
\begin{multline}\label{grothgp1} 
[H^1(\pi_1^t(\cau_s,\eta_s),M)]-[H^0(\pi_1^t(\cau_s,\eta_s),M)]\\ 
=(2g-2+\#(\cax_s\setminus\cau_s))[M]. 
\end{multline}
But the previous construction is functorial, so we can view
(\ref{grothgp1}) as an equality in the \textsc{Grothendieck} group of
cofinitely generated $\bbz_p[\cag]$-modules when
$M=F'_{p^{\infty}}\otimes_{\bbz}\bbz[\cag]$. Once again using
\textsc{Shapiro}'s lemma and lemma \ref{leminvpdiv} for the curve $\cau'$ we conclude that
$H^0(\pi_1^t(\cau_s,\eta_s),F'_{p^{\infty}}\otimes_{\bbz}\bbz[\cag])
=H^0(\pi_1^t(\cau'_s,\eta_s'),F_{p^{\infty}}')$ is finite, hence
$\bbz_p$-cotorsion, in particular it is killed by the functor $W$. It
then follows from (\ref{grothgp1}), (\ref{weq}) and (\ref{grothgp2})
that $W(S(\cac',F'_{p^{\infty}})$ is a $\bbq_p[\cag]$-free module of
rank equal to 
\begin{multline}\label{ineqtr}
(2d-2d_0')(2g-2+\#(\cax_s\setminus\cau_s))\\-\sum_{v\in\cac_s\setminus\cau_s}(2a_v+t_v)
+\#(\cac_s\setminus\cau_s)\cdot(2d_0')=\\
(2d-2d_0')(2g-2+\#(\cax_s\setminus\cac_s))\\+\sum_{v\in\cac_s\setminus\cau_s}(2a_v+t_v)
+f_{A,\cac}-\#(\cac_s\setminus\cau_s)\cdot(2d_0')\\
-\sum_{v\in\cac_s\setminus\cau_s}(2a_v+t_v)+\#(\cac_s\setminus\cau_s)(2d_0')\\
\le(2d-2d_0)(2g-2+\#(\cax_s\setminus\cac_s))+f_{A,\cac}.
\end{multline}
In order to see the validity of the latter inequality in (\ref{ineqtr}) we use the fact that the trace maps $\tau:B\to A$ and $\tau':B'\to A$ are injective, if $k$ has characteristic zero, and have finite kernel if $k$ has positive characteristic (for a discussion on this \emph{cf.} \cite[\S2 and \S6]{bc}). In particular, the induced map $B\to B'$ has necessarily finite kernel, therefore $d_0\le d'_0$ and this implies the searched inequality.
\end{proof}

Let $A':=A\times_KK'$, $\frF_{A'}$
the conductor divisor of $A'$ on $\cax'$ and
$f_{A'}:=\deg(\frF_{A'})$. Let 
$\frF_{A',\cac'}$ be the restriction of $\frF_{A'}$ to $\cac'$ and
$f_{A',\cac'}:=\deg(\frF_{A',\cac'})$. 

As observed in \cite{ell}, when $G$ is abelian and $H=\gal(l/k)$, then
$\epsilon(G,H)=\#\frO_{G_k}(G)$. Therefore, 
theorem \ref{mainthm} implies the following result.

\begin{theorem}\label{mainthm1} (cf. \cite[corollary 2.13]{ell})
With the same notation of theorem \ref{mainthm}, let $A/K$ be a non constant abelian variety of dimension $d$. Suppose that $k$ has either characteristic 0 or greater than $2d+1$. 
Let $\psi:\cac'\to\cac$ be a finite morphism defined over $k$ which is
geometrically abelian and \'etale with automorphism 
group $\cag=\aut_{k^s}(\cac'/\cac)$. Denote by $g'$ the genus of $\cac'$, $\cax'$ a regular compactification of $\cac'$ and $f_{A',\cac'}$ the conductor of $A'=A\times_KK'$ with respect to $\cac'$. Then
\begin{equation}\label{mainbd}
\rk\left(\frac{A(K')}{\tau'B'(k)}\right)\le
\frac{\#\frO_{G_k}(\cag)}{\#\cag}((2d-2d_0)(2g'-2+\#(\cax'_s\setminus\cac'_s))+f_{A',\cac'}).
\end{equation}
\end{theorem}

\begin{remark}\label{apell}
When $\cac=\cax$ is a complete curve and $k$ is a number field, under the
hypotheses of theorem \ref{mainthm}, (\ref{mainbd})
improves \cite[(1.7)]{pac}. Observe that here we make no hypothesis
concerning the truth of \textsc{Tate}'s conjecture and nor the irreducibility of
monodromy representations. Nevertheless, the method of \cite{pac}
allowed us to treat the
case of arbitrary ramification.

We now compare the bound obtained here with that of \cite[(1.4)]{pac}
in the ramified case. Let $k$ be a number field, $\cax/k$ be a
complete geometrically connected smooth curve defined over $k$ with
function field $K:=k(\cax)$. Let $\psi:\cax'\to\cax$ be a geometrically abelian cover defined
over $k$. Let $\car$ be the ramification locus of $\psi$,
$\cac:=\cax\setminus\car$, $\car':=\psi^{-1}(\car)$ and
$\cac':=\cax'\setminus\car'$. Then the restriction of $\psi$ to
$\cac'$ gives a geometrically abelian, now also geometrically \'etale, cover
$\cac'\to\cac$ of affine curves defined over $k$. 
Let $\car_s:=\car\times_{\spec(k)}\spec(k^s)=\cac_s\setminus\cau_s$. Let
$\frF_{A,\car}:=\frF_A-\frF_{A,\cac}$ and
$f_{A,\car}:=\deg(\frF_{A,\car})$. 

In
\cite[(1.4)]{pac}, under the aforementioned hypotheses, the bound obtained
was
\begin{equation}\label{oldbound}
\frac{\#\frO_{G_k}(\cag)}{\#\cag}(d(2d+1)(2g'-2))+\#\frO_{G_k}(\cag)\cdot
2d{f}_A.
\end{equation}
In particular, by \cite[proposition 3.7]{pac}, (\ref{oldbound}) is greater or equal to 
\begin{equation}\label{apchi}
\frac{\#\frO_{G_k}(\cag)}{\#\cag}(d(2d+1)(2g'-2)+2d{f}_{A'}).
\end{equation}
We now compare (\ref{apchi}) with (\ref{mainbd}). We see that
$d(2d+1)(2g'-2)\ge(2d-2d_0)(2g'-2)$, if $g'\ge1$. Clearly
$2df_{A',\cac'}\ge f_{A',\cac'}$. All we need to analyze is when
$2df_{A',\car'}\ge(2d-2d_0)\cdot\#\car'_s$. This inequality holds if
and only if for every $w\in\car'_s$ we have $2de_w\ge2d-2d_0$. The
latter inequality holds as long as $e_w\ge1$, \emph{i.e.}, $A'$ has
bad reduction at $w$. In fact, otherwise $d_0\ge d$, whence $d_0=d$. Since the characteristic of $k$ is zero, then the map $\tau:B\to A$ is injective, therefore $A$ would be constant. However we are ruling out
this possibility since the beginning. So theorem \ref{mainthm1} gives
a smaller bound for the rank of $A(K')/\tau'B'(k)$ than that of 
\cite[(1.7)]{pac} if 
\begin{enumerate}
\item[(a)] $g'\ge1$; and
\item[(b)] $\car'\subset\Delta'$, where $\Delta'$ denotes the discriminant
locus of $\phi':\caa'\to\cac'$. 
\end{enumerate}
\end{remark}

\section{Towers of function fields}

Let $\cac$ be a smooth geometrically connected curve defined over a field $k$ of characteristic $q$ which is either zero or greater than $2d+1$. We define a \emph{tower of curves} over $\cac$ to be a sequence 
$$
\cat : \cdots\to\cac_n\to\cdots\to\cac_1\to\cac_0:=\cac
$$
of finite geometrically \textsc{Galois} and \'etale covers of curves
$\cac_n\to\cac_0$ defined over $k$. For each cover $\cac_n\to\cac$
denote by $\cag_n:=\aut_{k^s}(\cac_n/\cac)$ the corresponding \textsc{Galois}
group. The \textsc{Galois} group of the tower $\cat$ is defined as
$\cag_{\infty}:=\varprojlim_n\cag_n$. 

Let $A/K$ be a non constant abelian variety. For each
$n\ge0$, let $K_n:=k(\cac_n)$ be the function field of $\cac_n$ and
$(\tau_n,B_n)$ be the $K_n/k$-trace of $A$. 

When we consider the question of the size of 
the rank of \textsc{Lang-N\'eron} groups of abelian
varieties over function fields we can ask the following 
vertical question : how does the rank of $A(K_n)/\tau_nB_n(k)$ vary
along the tower $\cat$? 

In the case where $\cac$ is a complete curve and $k$ is a number
field we studied this question in the special cases of the two
particular towers: 
\begin{enumerate}
\item[(*)] when $\cac$ is an elliptic curve and $\cac_n$ is obtained
  from $\cac$ as the pull-back by the multiplication by $n$ map in $\cac$;
\item[(**)] for a curve $\cac$ of any positive genus, 
$\cac_n$ is obtained as the pull-back of
  $\cac$ by the multiplication by $n$ map in the Jacobian variety
  $J_{\cac}$ of $\cac$. 
\end{enumerate}
Observe that the first situation was already dealt with by \textsc{Silverman}
in the case of elliptic fibrations in \cite{sil}. We proved in
\cite[theorems 6.2 and 6.5]{pac} that average rank of 
$A(K_n)/\tau_nB_n(k)$ as $n\to\infty$ was smaller than
a fixed multiple of the degree $f_A$ of the conductor of
$A$, under the hypothesis in (*) 
that $\cac$ had no complex multiplication and in (**) that its
Jacobian variety $J_{\cac}$ had $\ov{k}$-endomorphism ring equal to $\bbz$ (plus
additional hypotheses, \emph{cf.} \cite[theorem 6.5]{pac}).

Let $K_{\infty}:=\varinjlim_nK_n$ and $\tau_{\infty}B_{\infty}(k):=
\varinjlim_n\tau_nB_n(k)$. One may
naturally ask whether the abelian group 
$A(K_{\infty})/\tau_{\infty}B_{\infty}(k)$ is finitely generated. 

In \cite{ell}, \textsc{Ellenberg} considered this question in the case of a non constant
elliptic curve $E/K$ and supposed that each $\cag_n$ was isomorphic to
$\bbz/p^n\bbz$ so that $\cag_{\infty}\cong\bbz_p$. 
For each $n\ge0$, let
$\cak_n:=k^s(\cac_n)$ and $\cak_{\infty}:=\varinjlim_n\cak_n$. Then
under certain conditions \cite[theorem 4.4]{ell} he proved that
$E(\cak_{\infty})$ is finitely generated. The goal of this section is
to extend this result to the case of an abelian fibration. 

We fix a prime number $p$ different from $q$. We suppose that each $\cag_n$ is a finite $p$-group, whence $\cag_{\infty}$ is a pro-$p$ group. 
For every $n\ge0$, let
$F_{n,p^{\infty}}:=A[p^{\infty}]/\tau_nB_n[p^{\infty}]$, 
$\tau_{\infty}B_{\infty}[p^{\infty}]:=\varinjlim_n\tau_nB_n[p^{\infty}]$ and 
$F_{\infty,p^{\infty}}:=A[p^{\infty}]/\tau_{\infty}B_{\infty}[p^{\infty}]$. Let
$S(\cac_{\infty},F_{\infty,p^{\infty}}):=\varinjlim_nS(\cac_n,F_{n,p^{\infty}})$.
Then $S(\cac_{\infty},F_{\infty,p^{\infty}})$ is a discrete \linebreak 
$p$-primary group with a continuous action of $\cag_{\infty}$. Hence it also
comes equipped with an action of the \textsc{Iwasawa} algebra
$\Lambda(\cag_{\infty}):=\varprojlim_{\cah}\bbz_p[\cag_{\infty}/\cah]$, where $\cah$ runs through
the open normal subgroups of $\cag_{\infty}$. For every \'etale sheaf $\caf$
on $\cac$, given an integer $n\ge1$, denote by $\caf_{|\cac_n}$ the pull-back of $\caf$ to $\cac_n$ and let
$H^i_{\et}(\cac_{s,\infty},\caf):=\varinjlim_nH^i_{\et}(\cac_{s,n},\caf_{|\cac_n})$. 

\begin{hypothesis}\label{hypiw}
Assume that $\cag_{\infty}$ is a non trivial pro-$p$ finite dimensional
$p$-adic \textsc{Lie} group without $p$-torsion elements.
\end{hypothesis}

Under the hypothesis \ref{hypiw}, the \textsc{Iwasawa} algebra $\Lambda(\cag_{\infty})$ is
both right and left noetherian local ring without zero
divisors. Moreover, it follows from \cite[lemma 1.6]{how} that for every
cofinitely generated discrete $\Lambda(\cag_{\infty})$-module $M$ the \textsc{Galois} cohomology group
$H^i(\cag_{\infty},M)$ is a cofinitely generated
$\bbz_p$-module. Furthermore, it makes sense to define the
$\Lambda(\cag_{\infty})$-corank of $M$ as follows \cite{how}
$$
\corank_{\Lambda(\cag_{\infty})}(M):=\sum_{i\ge0}(-1)^i\corank_{\bbz_p}(H^i(\cag_{\infty},M)).
$$

\begin{lemma}
The $\Lambda(\cag_{\infty})$-module $S(\cac_{\infty},F_{\infty,p^{\infty}})$
is cofinitely generated and has $\Lambda(\cag_{\infty})$-corank equal to
$(2d-2d_0)(2g-2+\#(\cax_s\setminus\cac_s))+f_{A,\cac}$. 
\end{lemma}

\begin{proof}
The proof follows as in \cite[propositions 3.3 and 3.4]{ell} replacing
\cite[proposition 2.5]{ell} by proposition \ref{corkprop}. Observe
though that in the course of the proof of \cite[proposition 3.4]{ell}
it is necessary to have $E[p^{\infty}]^{\pi_1(\cau_s,\eta_s)}$
finite. In his case this followed from \cite[remark 2.4]{ell}, in the
current situation this follows from lemma \ref{leminvpdiv}.
\end{proof}

We assume from now on that $\cag_{\infty}\cong\bbz_p$. Let
$\cah_{\infty}:=\Ker(\pi_1^t(\cau_s,\eta_s)\twoheadrightarrow\cag_{\infty})$.

\begin{proposition}\label{fininvpdivpadic}
If $k$ is finitely generated over its prime
field $\kappa_0$, then the subspace 
$F_{0,p^{\infty}}^{\cah_{\infty}}$ of $F_{0,p^{\infty}}$
of the elements fixed under the action of $\cah_{\infty}$ is finite.
\end{proposition}

\begin{proof}
Assume first that $k$ is a finite field of characteristic $q$ and let
$\Frob_k$ be its \textsc{Frobenius} automorphism. Let $\pi_1^t(\cau,\eta_s)$ be
the \emph{arithmetic tame fundamental group} of $\cau$ with respect to
$\eta_s$. By the definition of the $\cah_{\infty}$ we have a short
exact sequence of groups
\begin{equation}\label{zpquot}
1\to\cah_{\infty}\to\pi_1^t(\cau_s,\eta_s)\to\cag_{\infty}\to1.
\end{equation} 
By the definition of $\pi_1^t(\cau,\eta_s)$, there is also another
sequence 
\begin{equation}\label{argeom}
1\to\pi_1^t(\cau_s,\eta_s)\to\pi_1^t(\cau,\eta_s)\to G_k\to1.
\end{equation} 
Since $\Frob_k$ acts on $\cah_{\infty}$, both sequences (\ref{zpquot})
and (\ref{argeom}) yield a third exact sequence
\begin{equation}\label{arquot}
1\to\cah_{\infty}\to\pi_1^t(\cau,\eta_s)\to\cag_{\infty}\ltimes G_k\to1.
\end{equation}

Suppose that $F_{0,p^{\infty}}^{\cah_{\infty}}$ is infinite. Let 
$V:=\Hom(F_{0,p^{\infty}}^{\cah_{\infty}},\bbq_p/\bbz_p)\otimes_{\bbz_p}\bbq_p$. This is a
finite positive dimensional $\bbq_p$-vector space endowed with an action of
$\cag_{\infty}\ltimes G_k$. The group $\cag_{\infty}$ acts on $V$
through its inclusion in $\cag_{\infty}\ltimes G_k$. Since
$\cag_{\infty}$ is abelian, then $V$ decomposes as a
$\bbq_p[\cag_{\infty}]$-module into a sum of one dimensional
eigenspaces. 

Let $\chi:\cag_{\infty}\to\bbq_p^*$ be a character of
$\cag_{\infty}$ and $V_{\chi}$ the eigensubspace of $V$ corresponding
to $\chi$. Let $V_{\chi}^{\Frob_k}$ be the subspace of $V_{\chi}$ of
those elements which are fixed by $\Frob_k$. Since the action of
$\sigma\in\cag_{\infty}$ on $V_{\chi}$ is through multiplication by
$\chi(\sigma)$, then $V_{\chi}^{\Frob_k}=V_{\chi^q}$. Similarly, for
every integer $n\ge1$, we have
$V_{\chi}^{\Frob_k^n}=V_{\chi^{q^n}}$. But $V$ is finite dimensional,
therefore there exists an integer $f\ge1$ such that $\chi^{q^f}=\chi$,
\emph{i.e.}, $\chi^{q^f-1}$ is the trivial character. However,
$\cag_{\infty}\cong\bbz_p$ is free, thus $\chi$ itself must be
trivial. In particular, $\cag_{\infty}$ acts trivially on $V$, whence
on $F_{0,p^{\infty}}^{\cah_{\infty}}$. In particular, the action of
$\cag_{\infty}\ltimes G_k$ on $F_{0,p^{\infty}}^{\cah_{\infty}}$
  reduces to that of $G_k$. Observe that by (\ref{argeom}) this
  only happens if and only if
$F_{0,p^{\infty}}^{\cah_{\infty}}
=(F_{0,p^{\infty}}^{\cah_{\infty}})^{\pi_1^t(\cau_s,\eta_s)}\subset 
F_{0,p^{\infty}}^{\pi_1^t(\cau_s,\eta_s)}$. 
However, the latter group is finite by lemma \ref{leminvpdiv}, and this yields a
contradiction. 

Suppose now that $k$ is a number field. For almost all prime ideals
$\frq$ of the ring of integers $\cao_k$ of $k$, the algebraic
varieties $\caa$ and $\cac$ reduce to smooth varieties $\caa_{\frq}$
and $\cac_{\frq}$ over the residue field $\bbf_{\frq}$ of
$\frq$. Moreover, the choice of $\frq$ can also be made so that
$\phi:\caa\to\cac$ reduces to a proper flat morphism
$\phi_{\frq}:\caa_{\frq}\to\cac_{\frq}$ also defined over
$\bbf_{\frq}$. Let $K_{\frq}:=\bbf_{\frq}(\cac_{\frq})$ be the
function field of $\cac_{\frq}$ and $K_{\frq}^s$ a separable closure
of $K_{\frq}$. The generic fiber of $\phi_{\frq}$
will be a non constant abelian variety $A_{\frq}$ defined over
$K_{\frq}$. 
Let $(\tau_{\frq},B_{\frq})$ be the $K_{\frq}/\bbf_{\frq}$-trace of $A_{\frq}$. 

By proper base change \cite[chapter VI, corollary
2.7]{mil} 
\begin{equation}\label{prbsch}
H^1_{\et}(A\times_{\spec(K)}\spec(K^s),\bbz_p)\cong 
H^1_{\et}(A_{\frq}\times_{\spec(K_{\frq})}\spec(K_{\frq}^s),\bbz_p).
\end{equation} 
It follows from \cite[theorem
15.1]{milav}
that 
\begin{equation}\label{pdiv0}
\begin{aligned}
&H^1_{\et}(A\times_{\spec(K)}\spec(K^s),\bbq_p/\bbz_p)\\
&\cong H^1_{\et}(A\times_{\spec(K)}\spec(K^s),\bbz_p)\otimes_{\bbz_p}\bbq_p/\bbz_p\\ 
&\cong\Hom_{\bbz_p}(T_p(A),\bbz_p)\otimes_{\bbz_p}\bbq_p/\bbz_p
\cong\Hom_{\bbq_p/\bbz_p}(A[p^{\infty}],\bbq_p/\bbz_p)\cong
A[p^{\infty}],
\end{aligned}
\end{equation}
where the latter isomorphism is not canonical, it is just an abstract 
isomorphism of $\bbq_p/\bbz_p$-modules. Similarly,
\begin{equation}\label{pdivq}
H^1_{\et}(A_{\frq}\times_{\spec(K_{\frq})}\spec(K_{\frq}^s),\bbq_p/\bbz_p)\cong 
A_{\frq}[p^{\infty}].
\end{equation}
It now follows from (\ref{prbsch}), (\ref{pdiv0}) and (\ref{pdivq}) that
$A[p^{\infty}]\cong A_{\frq}[p^{\infty}]$. The same argument used above 
also gives $\tau
B[p^{\infty}]\cong\tau_{\frq}B_{\frq}[p^{\infty}]$. \emph{A fortiori}, 
\begin{equation}\label{pdivcong}
F_{0,p^{\infty}}=A[p^{\infty}]/\tau B[p^{\infty}]\cong
A_{\frq}[p^{\infty}]/\tau_{\frq}B_{\frq}[p^{\infty}]=:F_{\frq,0,p^{\infty}}.
\end{equation}

Let $\cau_{\frq}$ be the reduction of $\cau$ modulo $\frq$ (it will be
equal to the smooth locus of $\phi_{\frq}$ for a generic choice of $\frq$) and $\eta_{\frq,s}$ a
geometric generic point of $\cau_{\frq}$. The reduction
$\cau_s\to\cau_{\frq,s}$ modulo $\frq$ induces the specialization
homomorphism
$\text{sp}_{\frq}:\pi_1^t(\cau_s,\eta_s)\cong\pi_1(\cau_s,\eta_s)
\to\pi_1^t(\cau_{\frq,s},\eta_{\frq,s})$ 
at the level of tame fundamental groups. This homomorphism is
necessarily surjective by \cite[exp. XIII, corollaire 2.12]{sga1}. 

Let $\cag_{\infty,\frq}$ be the Galois group of a geometric
$\bbz_p$-extension of the function field $\ov{\bbf}_{\frq}(\cac)$. Then we
have the following commutative diagram 
$$
\begin{CD}
1@>>>\cah_{\infty}@>>>\pi_1(\cau_s,\eta_s)@>>>\cag_{\infty}@>>>1\\
@.@V{\text{sp}'_{\frq,\infty}}VV@V{\text{sp}_{\frq}}VV@V{\text{sp}_{\frq,\infty}}VV\\
1@>>>\cah_{\frq,\infty}@>>>\pi_1^t(\cau_{\frq,s},\eta_{\frq,s})@>>>\cag_{\frq,\infty}@>>>1
\end{CD}.
$$

Observe that the the group isomorphism $\text{sp}_{\frq,\infty}$ is 
actually obtained via the specialization homomorphism
$\text{sp}_{\frq}$. In fact, $\text{sp}_{\frq}$ induces an isomorphism
$\text{sp}_{\frq}^{(q')}:\pi_1(\cau_s,\eta_s)^{(q')}
\overset{\cong}\longrightarrow\pi_1^t(\cau_{\frq,s},\eta_{\frq,s})^{(q')}$
between the maximal prime to $q$ quotients of both tame fundamental
groups (\emph{cf.} \cite[exp. X, corollaire 3.9]{sga1}). Consequently,
$\text{sp}_{\frq,\infty}$ is obtained from $\text{sp}^{(q')}_{\frq}$
by taking quotients on both sides. In particular, the diagram is
commutative. Finally, a simple diagram chasing then shows 
that the first vertical arrow is also a surjection. 

The action of $\cah_{\frq,\infty}$ on $F_{\frq,0,p^{\infty}}$ and the isomorphism
$F_{0,p^{\infty}}\cong F_{\frq,0,p^{\infty}}$, induce an action of
$\cah_{\frq,\infty}$ on $F_{0,p^{\infty}}$, thus
$\ker(\cah_{\infty}\twoheadrightarrow\cah_{\frq,\infty})$ acts
trivially on $F_{0,p^{\infty}}$. Therefore,
$F_{0,p^{\infty}}^{\cah_{\infty}}\cong
F_{0,p^{\infty}}^{\cah_{\frq,\infty}}\cong
F_{\frq,0,p^{\infty}}^{\cah_{\frq,\infty}}$. However, the latter group
is finite, by the first part of the proof. \emph{A fortiori}, the first
one is also finite. 

The same argument using the specialization homomorphism of tame
fundamental groups also works if $k$ is a one variable function field
over a finite field. In fact, one need only to notice that the surjectivity of
the specialization homomorphism of tame fundamental groups holds in
general. This follows from \cite[exp. XIII, corollaire 2.8]{sga1} and
it is enough to extend the result to the one variable function field case. 
Finally, an induction argument on the
transcendence degree of $k$ over its prime field $\kappa_0$, using the
specialization homomorphism, allows us to extend the result to any
field $k$ finitely generated over $\kappa_0$. 
\end{proof}

For each finite group $\cag_n$,  let $\cag^{(n)}:=\Ker(\cag_{\infty}\twoheadrightarrow\cag_n)$ and
$k_{\infty}$ the minimal algebraic extension of $k$ over which all
elements of $\cag_{\infty}$ are defined. As in \cite[proposition
4.1]{ell} we have the following result.

\begin{proposition}\label{propprofgp}
Suppose that $k$ is a finitely generated field over its prime field $\kappa_0$ of characteristic either zero or $q>2d+1$, where $d=\dim(A)$. Assume also that 
$p>(2d-2d_0)(2g-2+\#(\cax_s\setminus\cac_s))+f_{A,\cac}$ and $p\ne q$. Then there
exists an extension $l/k_{\infty}$ such that 
\begin{enumerate}
\item $\gal(k^s/l)$ acts trivially on
  $(A(k^s(\cac_n))/\tau_nB_n(k^s))\otimes_{\bbz}\bbq_p/\bbz_p$ for every 
  $n\ge0.$
\item $l$ is an abelian pro-$p$ extension of a finite extension of
  $k_{\infty}$. 
\end{enumerate}
\end{proposition}

\begin{proof}
The proof follows similarly as in \cite[proposition 4.1]{ell}, however
one has to consider the following point. Applying the
\textsc{Hochschild-Serre} 
spectral sequence to the tower $\cat$ one gets for every $n\ge0$ a
group homomorphism $f:S(\cac_n,F_{n,p^{\infty}})\to
S(\cac_{\infty},F_{\infty,p^{\infty}})^{\cag_{n}}$ whose kernel is
  equal to
$$
H^1(\cag^{(n)},H^0_{\et}(\cac_{\infty,s},\caf_{0.\infty}))
=H^1(\cag^{(n)},F_{0,p^{\infty}}^{\cah_{\infty}}).  
$$
But, by proposition
  \ref{fininvpdivpadic}, $F_{0,p^{\infty}}^{\cah_{\infty}}$ is
  finite. This replaces the argument of \cite[remark 2.4]{ell} in the
  proof of the proposition. 
\end{proof}

As a consequence of proposition \ref{propprofgp} we obtain 
 the following result.

\begin{theorem}\label{fginfty}(\emph{cf.} \cite[theorem 4.4]{ell}) Let $K=k(\cac)$ be the function field of a smooth geometrically connected curve. Let $A/K$ be a non constant abelian variety of dimension $d$. Assume that $k$ has characteristic either zero or $q>2d+1$. Suppose furthermore that
$p>(2d-2d_0)(2g-2+\#(\cax_s\setminus\cac_s))+f_{A,\cac}$ and $p\ne q$. 
Under the additional hypothesis : 
\begin{enumerate}
\item[($\dagger$)] for every extension $l/k$ which is an abelian pro-$p$ extension
of a finite extension of $k_{\infty}$, no divisible subgroup of
$S(\cac,F_{0,p^{\infty}})$ is fixed by $\gal(k^s/l)$; 
\end{enumerate}
the abelian group $A(k^s(\cac_{\infty}))/\tau_{\infty}B_{\infty}(k^s)$ is finitely generated. \end{theorem}

\begin{remark}
Once again it is necessary to replace \cite[remark 2.4]{ell} in the
proof of \cite[theorem 4.4]{ell} by proposition \ref{fininvpdivpadic} to get
theorem \ref{fginfty}.
\end{remark}

\begin{remark}\label{remgalmod}
In \cite[remark 4.5]{ell} it is discussed abstractly condition
($\dagger$). More precisely,  given a cofinitely generated
  $\bbz_p$-module $M$ with an action by $\gal(k^s/k_{\infty})$,
  condition ($\dagger$) means that the following property is satisfied. For
  every extension $l$ of $k_{\infty}$ which is an abelian pro-$p$
  extension of a finite extension of $k_{\infty}$, no divisible
  submodule of $M$ is fixed by $\gal(k^s/l)$. This property is
  inherited by submodules of $M$ as well as quotients of $M$ by finite
  submodules. It also respects exact sequences of modules $0\to M\to
  M'\to M''$ in the sense that if it holds for $M$ and $M''$, then it also
  holds for $M'$. 
\end{remark}

\section{Jacobian fibrations}
\subsection{Generalities}
In \cite[\S5]{ell} an example was given in which condition ($\dagger$)
is fulfilled for minimal elliptic K3 surfaces. It is natural to
ask whether such an example can be produced for higher dimensional
abelian fibrations. In this section we give a necessary condition for the existence of such an example in the context of Jacobian fibrations. However, due to the lack of examples of families of surfaces whose monodromy is ``sufficiently large'' (\emph{cf.} subsection \ref{famsurf}) we were not able to produce a concrete example as in the case of elliptic fibrations. 

In this section we
will assume that $\cac$ is a complete smooth geometrically irreducible curve defined over a subfield $k$ of
$\bbc$. For any variety $\cay$ defined over $k$, denote
$\ov{\cay}:=\cay\times_{\spec(k)}\spec(\ov{k})$. If $\caz$ is a
variety defined over $K:=k(\cac)$, denote
$\ov{\caz}:=\caz\times_{\spec(K)}\spec(\ov{K})$. 

Let $\cax$ be a smooth projective geometrically irreducible surface defined over $k$
and $\phi:\cax\to\cac$ be a proper flat morphism also defined over
$k$. The generic fiber $X$ of $\phi$ is a smooth projective
geometrically irreducible curve defined over $K$ of genus 
$d$ which we assume to be at least 2. The Jacobian fibration
$\phi_J$ associated to $\phi$ is a 
proper flat morphism $\phi_J:\caj\to\cac$ defined over $k$ from a
smooth geometrically irreducible $(d+1)$-fold 
$\caj$ defined over $k$ whose generic fiber is the Jacobian variety
$A:=\text{Jac}(X)$ of $X$. It has the property that for every $v\in\cac$ for which the fiber $\cax_v=\phi^{-1}(v)$ is  smooth, then the fiber $\phi_J^{-1}(v)$ coincides with the Jacobian variety $\text{Jac}(\cax_v)$ of $\cax_v$. Let $(\tau,B)$ be the $K/k$-trace of $A$. 

As before, given a prime number $p$, let
$F_{p^{\infty}}:=A[p^{\infty}]/\tau B[p^{\infty}]$. Let $\ov{\eta}$ be the
geometric generic point of $\cac$, $j:\ov{\eta}\hookrightarrow\ov{\cac}$
the inclusion map and $\caf_{p^{\infty}}:=j_*(F_{p^{\infty}})$. Let
$\tilde{\caf}_{p^{\infty}}:=R^1\phi_*(\bbq_p/\bbz_p)$. 

Suppose that $(\tau,B)$ is trivial. It follows from \cite[corollary
9.6]{miljv} after tensoring with $\bbq_p/\bbz_p$ that
$$
H^1_{\et}(\ov{X},\bbq_p/\bbz_p)\cong 
H^1_{\et}(\ov{A},\bbq_p/\bbz_p). 
$$ 
It then follows from
\cite[theorem 15.1]{milav}
that 
\begin{multline*}
H^1_{\et}(\ov{A},\bbq_p/\bbz_p)\cong H^1_{\et}(\ov{A},\bbz_p)\otimes_{\bbz_p}\bbq_p/\bbz_p
\cong\Hom_{\bbz_p}(T_p(A),\bbz_p)\otimes_{\bbz_p}\bbq_p/\bbz_p\\
\cong \Hom_{\bbq_p/\bbz_p}(A[p^{\infty}],\bbq_p/\bbz_p)\cong
A[p^{\infty}], 
\end{multline*}
where the latter isomorphism is not canonical, just as
abstract $\bbq_p/\bbz_p$-modules. As a consequence,
\begin{equation}\label{genisosh}
j^*(\tilde{\caf}_{p^{\infty}})=(\tilde{\caf}_{p^{\infty}})_{\ov{\eta}}
=H^1_{\et}(\ov{X},\bbq_p/\bbz_p)\cong
H^1_{\et}(\ov{A},\bbq_p/\bbz_p)\cong A[p^{\infty}]=F_{p^{\infty}},
\end{equation}
whence
$j_*(j^*(\tilde{\caf}_{p^{\infty}}))\cong\caf_{p^{\infty}}$.
Therefore, we have a surjective map $H^1_{\et}(\ov{\cac},\tilde{\caf}_{p^{\infty}})\to
H^1_{\et}(\ov{\cac},\caf_{p^{\infty}})=S(\cac,F_{p^{\infty}})$, since the kernel of 
$\tilde{\caf}_{p^{\infty}}\to\caf_{p^{\infty}}$ has zero dimensional support.

The \textsc{Leray} spectral sequence
$H^i_{\et}(\ov{\cac},R^j\phi_*(\bbq_p/\bbz_p))\implies
H^{i+j}_{\et}(\ov{\cax},\bbq_p/\bbz_p)$ yields the exact sequence of
cohomology groups 
\begin{multline}\label{leray1}
0\to H^1_{\et}(\ov{\cac},\bbq_p/\bbz_p)\to
H^1_{\et}(\ov{\cax},\bbq_p/\bbz_p)\to
H^0_{\et}(\ov{\cac},\tilde{\caf}_{p^{\infty}})\to \\
H^2_{\et}(\ov{\cac},\bbq_p/\bbz_p)\overset{d_2}\longrightarrow 
H^2_{\et}(\ov{\cax},\bbq_p/\bbz_p).
\end{multline}
Let $F$ be a smooth fiber of $\phi$. Then the image of $d_2$ in
(\ref{leray1}) is generated by the class $[F]$ of $F$ in
$H^2_{\et}(\ov{\cax},\bbq_p/\bbz_p)$. Let $M$ be the quotient of the
latter group by the subgroup generated by $[F]$. Now the previous
spectral sequence at degree 2 yields
\begin{equation}\label{leray2}
0\to H^1_{\et}(\ov{\cac},\tilde{\caf}_{p^{\infty}})\overset{d_1}\longrightarrow M\to
H^0_{\et}(\ov{\cac},R^2\phi_*(\bbq_p/\bbz_p))\to H^2_{\et}(\ov{\cac},\tilde{\caf}_{p^{\infty}}).
\end{equation}
Observe that $H^2_{\et}(\ov{\cac},\tilde{\caf}_{p^{\infty}})$ is finite,
since it is dual to
$$
H^0_{\et}(\ov{\cac},\tilde{\caf}_{p^{\infty}})\cong
H^1_{\et}(\ov{X},\bbq_p/\bbz_p)^{\pi_1^t(\ov{\cau},\ov{\eta})}\cong
A[p^{\infty}]^{\pi_1^t(\ov{\cau},\ov{\eta})}, 
$$ 
where the latter
isomorphism follows from (\ref{genisosh}) and the finiteness of the
latter space is a consequence of lemma \ref{leminvpdiv}. Hence, 
\begin{equation}\label{corkforh2}
\corank_{\bbz_p}(H^1_{\et}(\ov{\cac},\tilde{\caf}_{p^{\infty}}))=\corank_{\bbz_p}(M)
-\corank_{\bbz_p}(H^0_{\et}(\ov{\cac},R^2\phi_*(\bbq_p/\bbz_p))).
\end{equation} 
Observe that the generic stalk of $R^2\phi_*(\bbq_p/\bbz_p)$ has corank
1. For each $v\in\ov{\cac}$ let $m_v$ be the number of irreducible
components of the fiber $\cax_v:=\phi^{-1}(v)$. Then the corank of
$R^2\phi_*(\bbq_p/\bbz_p)$ at $v$ equals $m_v-1$. Therefore
\begin{equation}\label{corkr2}
\corank_{\bbz_p}(H^0_{\et}(\ov{\cac},R^2\phi_*(\bbq_p/\bbz_p)))=1+\sum_{v\in\ov{\cac}}(m_v-1).
\end{equation}

Let $\pic(\cax)$ be the \textsc{Picard} group of $\ov{k}$-divisor classes of
$\cax$. The composition 
\begin{equation}\label{degree1}
M\to H^0_{\et}(\ov{\cac},R^2\phi_*(\bbq_p/\bbz_p))\to
H^0_{\et}(\ov{\cac},j_*(j^*(R^2\phi_*(\bbq_p/\bbz_p))))\cong\bbq_p/\bbz_p
\end{equation}
is the degree map. It sends the image in
$H^2_{\et}(\ov{\cax},\bbq_p/\bbz_p)$ of a class $[D]$ of a divisor in
$\pic(\cax)$ to its intersection number $(D\cdot X)$ with the generic
fiber. Let $G(\cax)$ be the quotient of the subset of elements of degree zero in
$H^2_{\et}(\ov{\cax},\bbq_p/\bbz_p)$ by the module generated by the
class of $[F]$. Thus,
$H^1_{\et}(\ov{\cac},\tilde{\caf}_{p^{\infty}})$ is a submodule of
$G(\cax)$. 

Let $(D_0\cdot X)$ be a generator of the ideal $\{(D\cdot
X)\,|\,D\in\Div(\cax)\}$, where $(\ \cdot\ )$ denotes the intersection
pairing on $\pic(\cax)$. Let $\imath:X\to\cax$ be the inclusion map
and $\imath^*:\Div(\cax)\to\Div(X)$ the pull-back map obtained by
restricting the divisors to the generic fiber $X$. Define
$\psi:\pic(\cax)\to A$ by 
$$
\psi([D]):=\imath^*\left(D-\frac{(D\cdot X)}{(D_0\cdot X)}D_0\right). 
$$
It follows from \cite[lemme 3.7]{hp} that $\psi$ is surjective and its
kernel $\cas$ is generated by the classes $[D_0]$, $[F]$ and the
classes of all irreducible components of the singular fibers of $\phi$
(except one). \emph{A fortiori}, 
\begin{equation}\label{rkker}
\rk_{\bbz}(\cas)=2+\sum_{v\in\ov{\cac}}(m_v-1).
\end{equation}
Let $\ns(\cax)$ be the \textsc{N\'eron-Severi} group of the surface $\cax$. 
We also extended the \textsc{Shioda-Tate} formula (\emph{cf.} \cite[proposition
3.8]{hp}) to fibrations not necessarily having a section and without
any hypothesis on the $K/k$-trace $(\tau,B)$ of $A$ being trivial. As
a consequence of this result we obtained 
\begin{equation}\label{shiota}
\rk_{\bbz}(\ns(\cax))
=\rk_{\bbz}\left(\frac{A(\ov{k}(\cac))}{\tau B(\ov{k})}\right)
+\rk_{\bbz}(\cas).
\end{equation}

The cohomology group $H^2_{\et}(\ov{\cax},\bbz_p)$ comes
equipped with a quadratic form $q_{\cax}$ through the pairing defined
by the cup product. Note also that for the images of the elements of 
$\pic(\cax)$ in $H^2_{\et}(\ov{\cax},\bbz_p)$, the intersection
pairing in $\pic(\cax)$ is compatible with the cup product in
$H^2_{\et}(\ov{\cax},\bbz_p)$. Moreover,
$\corank_{\bbz_p}(H^2_{\et}(\ov{\cax},\bbz_p))$ equals 
the second Betti number $b_2(\cax)$ of $\cax$. 
Let $\Gamma'_p$ be the subgroup of automorphisms of $H^2_{\et}(\ov{\cax},\bbz_p)$ which
preserves $q_{\cax}$ and stabilizes $[F]$ and $[D_0]$. Let $\Gamma_p$ be a
finite index subgroup of $\Gamma'_p$.  

\subsection{A sufficient condition}
\begin{theorem}\label{thmjac}(\emph{cf.} \cite[theorem 5.1]{ell}) 
Let $\cax$, resp. $\cac$, be a smooth projective irreducible surface,
resp. curve, defined over a number field $k$. 
Let $\phi:\cax\to\cac$ be a proper flat morphism also defined over
$k$. Let $d\ge2$ be the genus of
the generic fiber $X$ of $\phi$. Let $p$ be a prime number, $n\ge0$
an integer, $h$ a function in
$K:=k(\cac)$, $K_n:=K(h^{1/p^n})$, 
$\cak:=\ov{k}(\cac)$ and $\ov{K}_n:=\cak(h^{1/p^n})$. 
Let $\phi_J:\caj\to\cac$ be the
Jacobian fibration associated to $\phi$ and $A$ the generic fiber
of $\phi_J$. Suppose the following 
\begin{enumerate}
\item[(a)] The image of $\gal(\ov{k}/k)$ in
  $\aut(H^2_{\et}(\ov{\cax},\bbz_p))$ contains $\Gamma_p$.
\item[(b)] $d$ is either 2, 6 or odd and for every $v\in\ov{\cac}$
  which is not a zero nor a pole of $h$, the fiber $\cax_v$ of $\phi$ is smooth
  and $\End_{\ov{\kappa}_v}(\mathrm{Jac}(\cax_v))=\bbz$, where $\ov{\kappa}_v$ denotes an algebraic closure of the residue field $\kappa_v$ of $v$ and $\mathrm{Jac}(\cax_v)$ is the Jacobian variety of $\cax_v$. 
\item[(c)] $p>b_2(\cax)-2+2d\cdot\deg(h)$. 
\end{enumerate}
Then the rank of the \textsc{Mordell-Weil} group $A(\ov{K}_n)$ is uniformly 
bounded as $n\to\infty$. 
\end{theorem}

\begin{proof}We start by observing that condition (a) implies that for every $n\ge1$, the $K_n/k$-trace of $A$ is trivial. Let $k_{\infty}:=k(\zeta_{p^{\infty}})$ be the field obtained from $k$ by adjoining all $p$-th power roots of unity. By the previous condition, the image of $\gal(\ov{k}/k_{\infty})$ in $\aut(H^2_{\et}(\ov{\cax}),\bbz_p)$ still contains a finite index subgroup of $\Gamma_p$, since the determinant map sends $\Gamma_p$ to a finite group. Next, denote by $H$ the $\bbz_p$-module generated by the images of the classes of $F$ and $D_0$. Then $\gal(\ov{k}/k_{\infty})$ acts irreducibly on $(H^2_{\et}(\ov{\cax},\bbz_p)/H)\otimes_{\bbz_p}\bbq_p$. It follows from the second exact sequence in p.29 of \cite{ray} that we have a surjective map $\psi:H^2_{\et}(\ov{\cax},\bbz_p)\twoheadrightarrow H^2_{\et}(\ov{\cax},\bbg_m)\{p^{\infty}\}$, where the latter denotes the $p$-primary subgroup of $H^2_{\et}(\ov{\cax},\bbg_m)$. As a consequence $\gal(\ov{k}/k_{\infty})$ also acts irreducibly on $(H^2_{\et}(\ov{\cax},\bbg_m)\{p^{\infty}\}/\psi(H))\otimes_{\bbz_p}\bbq_p$. Recall the inclusion map $j:\ov{\eta}\hookrightarrow\cac$. It follows from the last map of p.28 of \cite{ray} that there exists an isomorphism $\vartheta:H^1_{\et}(\ov{\cac},j_*A)\{p^{\infty}\}\otimes_{\bbz_p}\bbq_p\overset{\cong}\longrightarrow H^2_{\et}(\ov{\cax},\bbg_m)\{p^{\infty}\}\otimes_{\bbz_p}\bbq_p$, where in the left hand side we are considering $A$ as an \'etale sheaf on $\ov{\eta}$. \emph{A fortiori}, $\gal(\ov{k}/k_{\infty})$ also acts irreducibly on $(H^1_{\et}(\ov{\cac},j_*A)\{p^{\infty}\}\otimes_{\bbz_p}\bbq_p)/\vartheta^{-1}(\psi(H)\otimes_{\bbz_p}\bbq_p)$. However, if the $K/k$-trace $(\tau,B)$ of $A$ is not zero, then this cannot happen, since the latter $\bbq_p[\gal(\ov{k}/k_{\infty})]$-module admits $(H^1_{\et}(\ov{\cac},j_*(\tau(B)))\{p^{\infty}\}\otimes_{\bbz_p}\bbq_p)/((H^1_{\et}(\ov{\cac},j_*(\tau(B)))\{p^{\infty}\}\otimes_{\bbz_p}\bbq_p)\cap\vartheta^{-1}(\psi(H)\otimes_{\bbz_p}\bbq_p)$ as a $\bbq_p[\gal(\ov{k}/k_{\infty})]$-submodule. Hence, $B=0$. The same argument applies the $K_n/k$-trace of $A$ for every $n\ge1$.

As a consequence of (\ref{corkforh2}) and (\ref{corkr2}) we have 
\begin{equation}\label{corkr1}\corank_{\bbz_p}(H^1_{\et}(\ov{\cac},\tilde{\caf}_{p^{\infty}}))=b_2(\cax)-2-\sum_{v\in\ov{\cac}}(m_v-1).\end{equation}
The second thing we need to do is to show that this corank equals $2d(2g-2)+f_A$. Indeed, it follows from the discussion on \cite[\S3]{ray} and the formula \cite[th\'eor\`eme 3, (ii)]{ray} that  
\begin{equation}\label{rayeq}b_2(\cax)-\rk_{\bbz}(\ns(\cax))=-\rk(A(\ov{k}(\cac)))+2d(2g-2)+f_A.\end{equation}
Hence, by (\ref{shiota}), (\ref{rkker}), (\ref{rayeq}) and \ref{corkr1}, we have 
\begin{equation}\label{eqirred}\corank_{\bbz_p}(H^1_{\et}(\ov{\cac},\tilde{\caf}_{p^{\infty}}))=2d(2g-2)+f_A.\end{equation}
(In fact the previous calculations would also hold if the $B\ne0$, however we would need to add $4\dim(B)$ to the right hand side of (\ref{eqirred}).) 

The rest of the proof is very similar to that of \cite[theorem 5.1]{ell}, we restrict ourselves to just pointing out the differences between them. As in \emph{loc. cit.}, $G(\cax)$ satisfies property ($\dagger$), hence  so does its submodule $H^1_{\et}(\ov{\cac},\tilde{\caf}_{p^{\infty}})$ (\emph{cf.} remark \ref{remgalmod}). Since the corank of this module equals $2d(2g-2)+f_A$ which is the corank of $S(\cac,F_{p^{\infty}})$, we conclude that the surjection
$H^1_{\et}(\ov{\cac},\tilde{\caf}_{p^{\infty}})\twoheadrightarrow 
H^1_{\et}(\ov{\cac},\caf_{p^{\infty}})$ has finite kernel. \emph{A
  fortiori} (\emph{cf.} remark \ref{remgalmod}),
$H^1_{\et}(\ov{\cac},\caf_{p^{\infty}})$ also satisfies condition
($\dagger$). 

Let $\ov{\caz}\subset\ov{\cac}$ be the scheme of zeroes and poles of
$h$ in $\ov{\cac}$. Let $\ov{\cav}:=\ov{\cac}\setminus\ov{\caz}$. 
Whence we have an exact sequence (\emph{cf.} \cite[chapter III,
proposition 1.25]{mil}) 
\begin{multline}\label{exschz} 
0\to H^1_{\et}(\ov{\cac},\caf_{p^{\infty}})\to
H^1_{\et}(\ov{\cav},\caf_{p^{\infty}})\to \\ 
H^0_{\et}(\ov{\caz},\caf_{p^{\infty}}(-1)_{|\ov{\caz}})=
\bigoplus_{v\in\ov{\caz}}H^1_{\et}(\cax_v,(\bbq_p/\bbz_p)(-1)),
\end{multline} 
where in the last equality of (\ref{exschz}) we used Poincar\'e's
duality for \'etale cohomology. 

It follows from hypothesis (a) that the action of $\gal(\ov{k}/k_{\infty})$ on
$G(\cax)$ is irreducible. Hence all fibers of $\phi$ are irreducible.
It follows from \cite[proposition 9.6]{miljv} 
that there exists an isomorphism (after having tensored both sides by $(\bbq_p/\bbz_p)(-1)$) 
$$
H^1_{\et}(\ov{\cax}_v,(\bbq_p/\bbz_p)(-1))\cong 
H^1_{\et}(J_{\ov{\cax}_v},(\bbq_p/\bbz_p)(-1)). 
$$
But by hypothesis (b), 
the \textsc{Mum\-ford-Tate} conjecture is true for $\text{Jac}(\ov{\cax}_v)$ (\emph{cf.}
remark \ref{remmt}). Consequently,
$H^0_{\et}(\ov{\caz},\caf_{p^{\infty}}(-1)_{|\ov{\caz}})$ also
satisfies property ($\dagger$). In particular, by remark \ref{remgalmod}, 
$H^1_{\et}(\ov{\cav},\caf_{p^{\infty}})=S(\cav,F_{p^{\infty}})$ also
satisfies property ($\dagger$). Moreover by (\ref{exschz}) its $\bbz_p$-corank is at
most equal to $b_2(\cax)-2+2d\cdot\deg(h)<p$. The choice of $\cav$ implies that the
map $h\mapsto h^{1/p^n}$ at the level of functions gives an \'etale Galois
cover $\cav_n\to\cav$ at the level of curves and we are exactly in the
set-up of theorem \ref{fginfty}, where $\cav$ plays the role of the
affine curve in that statement. The theorem is now a consequence of
the latter result.
\end{proof}

\begin{remark}\label{remmt}
We introduced condition (b) in the hypotheses of the theorem to assure the truth of the \textsc{Mumford-Tate} 
conjecture for the fibers $\caj_v$ of $\phi_J$ for $v\in\ov{\cac}$ outside the
support of the divisor of $h$. For a discussion on the
\textsc{Mumford-Tate} 
conjecture see \cite{pink}. For us here the only important thing is
the following consequence. Suppose that $Y$ is an abelian variety of
dimension $d$ 
defined over a number field $k$. Assume that $k$ is sufficiently large
so that $Y$ is principally polarized. For each integer $n\ge1$, let
$Y[n]$ be the subgroup of $n$-torsion points of $Y$. \textsc{Serre}
proved in \cite[th\'eor\`eme 3]{sercf} that if $\End_{\ov{k}}(Y)=\bbz$ and its dimension is either 2, 6 or odd, then
the image of the \textsc{Galois} representation 
$\rho_n:\gal(\ov{k}/k)\to\aut(Y[n])\cong\gsp_{2d}(\bbz/n\bbz)$ 
has index at most $I(Y,k)$
(independent of $n$) for every $n\ge1$ (\emph{cf.} \cite[theorem
6.3]{pac}). 
\end{remark}

\subsection{Towards families of surfaces}\label{famsurf}The hardest condition on theorem \ref{thmjac} is (a). The idea to produce an example where it might be satisfied is the following. Suppose that there exists a proper flat family of surfaces $\pi:\frX\to\frS$ parametrized by a smooth projective variety $\frS$. Assume that $\frX$, $\frS$ and $\pi$ are defined over $k$. Suppose furthermore that for every $s\in\frS$ the fiber $\frX_s$ of $\pi$ at $s$ comes equipped with a fibration on curves $\phi_s:\frX_s\to\cac_s$ to a smooth projective geometrically connected curve $\cac_s$. Let $L=k(\frS)$ be the function field of $\frS$, $\ov{L}$ an algebraic closure of $L$, $X/L$ the generic fiber of $\pi$ (which is also assumed to be projective and geometrically connected), $\ov{X}:=X\times_{\spec(L)}\spec(L)$ and $\eta$ a geometric generic point of $\frS$. Denote by $\pi_1^{\text{geom}}:=\pi_1(\ov{X},\eta)$ the geometric algebraic fundamental group of $X$ with respect to $\eta$. Consider the monodromy representation $\rho:\pi_1^{\text{geom}}\to\aut(H^2_{\et}(\ov{X},\bbz_p))$ and let $G^{\text{geom}}$ be its geometric monodromy group, \emph{i.e.}, the \textsc{Zariski} closure of the image of $\rho$ in the algebraic group $\gl_{N,\bbz_p}$, where $N=\rk_{\bbz_p}(H^2_{\et}(\ov{X},\bbz_p))$. If we can prove that $G^{\text{geom}}$ is the orthogonal algebraic group $O_{N,\bbz_p}$, then a similar argument to that in \cite[\S5]{ell} (supposing that $k$ is a number field and using the \textsc{Hilbert}'s irreducibility theorem) allows one to obtain for $u$ in an open subset $U$ of $\frS$ a representation $\rho_u:\gal(\ov{k}/k)\to\aut(H^2_{\et}(\ov{\frX}_u,\bbz_p))$ induced by $\rho$ with the following property (here $\frX_u$ denotes the element in the family $\pi$ corresponding to the fiber at $u$). The image of $\rho_u$ must contain a subgroup of finite index of the group of automorphisms of $H^2_{\et}(\ov{\frX}_u,\bbz_p)$ which preserves the quadratic form defined by the cup product and stabilizes a fiber $F$ of $\phi_u$ and horizontal divisor $D_0$ in $\frX_u$ (with notation as in the previous subsection). This would give an infinite number of surfaces satisfying condition (a). However, examples of families of surfaces whose geometric monodromy group is the orthogonal group seem still to be lacking.



\begin{thebibliography}{BoLuRa90}
\bibitem[BoLuRa90]{bolura}S. Bosch, W. L\"utkebohmert, M. Raynaud,
  \emph{N\'eron Models}, Springer-Verlag, 1990.
\bibitem[Co05]{bc}B. Conrad, \emph{The $K/k$-image, the $K/k$-trace of
    abelian varieties and a theorem of N\'eron and Lang}, 
    Eiseng. Math. \textbf{52} (2006), 37-108.
\bibitem[Del81]{del}P. Deligne, \emph{Conjectures de Weil II}, Pub. Math. IHES \textbf{52} (1981), 313-428.
\bibitem[De81]{mds}M. Deschamps, \emph{R\'eduction semi-stable}, in
  S\'em. sur les pinceux des courbes de genre au moins deux,
  Ast\'erisque 86 (1981), 1-34.
\bibitem[El06]{ell}J. Ellenberg, \emph{Selmer groups and Mordell-Weil groups
    of elliptic curves over
    towers of function fields}, Compositio Math. \textbf{142} (2006),
    1215-1230.
\bibitem[Gr71]{sga1}A. Grothendieck, \emph{S\'em. G\'eom\'etrie
    Alg\'ebrique 1}, Documents Math. \textbf{3}, Soc. Math. France, 2003.
\bibitem[Gr72]{sga7}A. Grothendieck, \emph{S\'em. G\'eom\'etrie Alg\'ebrique
  7, vol. I}, Lec. Notes Math. 288, Springer-Verlag, 1972.
\bibitem[HiPa05]{hp}M. Hindry, A. Pacheco, \emph{Sur le rang des
    Jacobiennes sur un corps de fonctions}, Bull. Soc. Math. France
    \textbf{133} (2005), 275-295.
\bibitem[Ho02]{how}S. Howson, \emph{Euler characteristics as invariants of
    Iwasawa modules}, Proc. London Math. Soc. \textbf{85} (2002),
    634-658.
\bibitem[Ka81]{katz}N. Katz, \emph{Monodromy of families of curves},
  in S\'em. Th. Nombres Paris 1979-80, pp. 64-84, 
Birkh\"auser, 1981.
\bibitem[Ka72]{katzsga7}N. Katz, \emph{Pinceax de Lefschetz:
    th\'eor\`eme d'existence}, SGA \textbf{7} II, Lect. Notes in
    Math. \textbf{340}, exp. VII, 212-253, Springer-Verlag, 1972.
\bibitem[Ka93]{katzjams}N. Katz, \emph{Affine cohomological transforms,
    perversity and monodromy}, J. Amer. Math. Soc. \textbf{6} (1993), 149-222.
\bibitem[La83]{la}S. Lang, \emph{Fundamentals of Diophantine
    Geometry}, Springer-Verlag, 1983.
\bibitem[Ma55]{mat}A. Mattuck, \emph{Abelian varieties over $p$-adic ground
    fields}, Annals of Math. \textbf{62} (1955), 92-119.
\bibitem[Mi80]{mil}J. S. Milne, \emph{\'Etale Cohomology}, Princeton University Press,
 Princeton, 1980.
\bibitem[Mi85jv]{miljv}J. S. Milne, \emph{Jacobian varieties} in
  \emph{Arithmetic Geometry}, eds. G. Cornell, J. H. Silverman,
  pp. 167-212, 1985.
\bibitem[Mi85av]{milav}J. S. Milne, \emph{Abelian varieties} in
  \emph{Arithmetic Geometry}, eds. G. Cornell, J. H. Silverman,
  pp. 103-150, 1985.
\bibitem[NeScWi00]{neucoh}J. Neukirch, A. Schmidt, K. Wingberg,
  \emph{Cohomology of Number Fields}, Springer-Verlag, Grundlehren, 2000.
\bibitem[Og62]{ogg}A. P. Ogg, \emph{Cohomology of abelian
    varieties over function fields}, Annals of Math. \textbf{76}
    (1962), 185-212.
\bibitem[Pa05]{pac}A. Pacheco, \emph{The rank of abelian varieties over
    function fields}, 2005, Manuscripta Math. \textbf{118} (2005),
    361-381.
\bibitem[Pi98]{pink}R. Pink, \emph{$\ell$-adic monodromy groups,
    cocharacters and the Mumford-Tate conjecture}, J. reine und
    angewandet Mathematik (Crelle) \textbf{495}, 187-237.
\bibitem[Ra68]{ray}M. Raynaud, \emph{Caract\'eristique
    d'Euler-Poincar\'e d'un faisceau et cohomologie des vari\'et\'es
    ab\'eliennes, (d'apr\`es Ogg-Shafarevich et Grothendieck)},
    S\'em. Bourbaki \textbf{279} (1964/65), in \emph{Dix expos\'es sur
    la cohomologie des sch\'emas}, eds. A. Grothendieck, et al.
\bibitem[Se86]{sercoh}J.-P. Serre, \emph{Cohomologie Galoisienne},
  Lecture Notes in Math. \textbf{5}, Springer-Verlag, 1986.
\bibitem[SeCF85]{sercf}J.-P. Serre, \emph{R\'esum\'e des cours},
  Coll\`ege de France, 1985/86, \OE uvres, Collected Papers IV,
  pp. 33-37.
\bibitem[SeTa68]{seta}J.-P. Serre, J. Tate, \emph{Good reduction of abelian
    varieties} Annals of Math. \textbf{88} (1968), 492-517.
\bibitem[Si02]{sil}J. Silverman, \emph{The rank of elliptic surfaces in unramified abelian towers over number fields},
J. reine und angewandet Mathematik (Crelle) \textbf{577} (2004),
153-169.
\bibitem[Sh99]{shi}T. Shioda, \emph{Mordell-Weil lattices for higher genus
fibration over a curve}, New Trends in Algebraic Geometry, London
Math. Soc. Lecture Notes \textbf{264}
(1999), 359-373.
\bibitem[Ta65]{tate}J. Tate, \emph{Algebraic cycles and poles of
zeta-functions}, Arithmetical Algebraic Geometry, Harper and Row, New York, 93-110 (1965).
\end{thebibliography}
\end{document}